
\font\elevenrm=cmr10 scaled\magstephalf
\font\ninerm=cmr9
\elevenrm
\hsize=30truecc
\vsize=44truecc
\raggedbottom
\medskip \vskip 4pt plus 1pt minus 1pt


\def\c{\cdot}
\def\ot{\otimes}
\def\rt{\rightarrow}
\def\ld{\ldots}

\def\ba{\begin{array}}
\def\ea{\end{array}}

\def\R{\hat{R}}
\def\I{\hat{I}}
\def\x{\hat{x}}
\def\A{\hat{A}}
\def\pih{\hat{\pi}}
\def\D{\hat{D}}
\def\d{\hat{d}}
\def\M{\cal{M}}
\def\al{\alpha}
\def\la{\lambda}
\def\ga{\gamma}
\def\de{\delta}
\def\be{\beta}
\def\1{x^1}
\def\2{x^2}


\def\title#1{\centerline{\bf #1}}
\def\author#1{\medskip\centerline{#1}}
\def\address#1{\centerline{\ninerm \it #1}}

\title{FIRST ORDER OPTIMUM CALCULI\footnote{}{\ninerm
{\it 1991 Mathematics Subject Classification: 16D20; 16W25;
16U80.}\hfill\break
This work has been supported by grants from State Research Committee
KBN
No 2 P302 02307 (A. B.) and Russian Foundation for Fundamental
Research
No 93-01-16171 (V. K. K.).}}
\author{ANDRZEJ BOROWIEC}
\address{Institute of Theoretical Physics, University of Wroc{\l}aw}
\address{Plac Maxa Borna 9, 50204 Wroc{\l}aw, Poland}
\address{E-mail: borowiec@ift.uni.wroc.pl}
\author{VLADISLAV K. KHARCHENKO}
\address{Institute of Mathematics, Novosibirsk, Russia}
\address{E-mail: kharchen@math.nsk.su}

\medskip{\ninerm {\bf Abstract:}\
A new notion of an optimum first order calculi was
introduced in [Borowiec, Kharchenko and Oziewicz, 1993].
A module of vector fields for a coordinate differential is defined.
Some examples of
optimal algebras for homogeneous bimodule commutations are presented.
Classification theorem for homogeneous calculi
with  commutative optimal algebras in two variables is
proved.}\medskip


\baselineskip13pt

{\bf 1. Introduction. }\
Quantum spaces are identified with noncommutative algebras.
Differential calculi on quantum spaces have been
elaborated by Pusz and Woro\-no\-wicz [6], by Pusz [5] and by
Wess and Zumino [7].
A bicovariant differential calculus on  quantum groups was presented
by
Woronowicz [8]. Woronowicz found that a
construction of a first order calculus, a bimodule of one forms,
is not functorial. There are many nonequivalent calculi for a given
associative algebra.
This yields a problem of classification of first order calculi.

In our previous paper [1] following the developments
by Pusz, Woro\-no\-wicz [6] and by Wess, Zumino [7],
we have proposed a general algebraic formalism for  first order
calculi on an arbitrary associative algebra with a given
presentation. A basic idea was that a noncommutative differential
calculus is best
handled by means of commutation relations between generators of an
algebra and its
differentials. We are assuming that a bimodule of one
forms is a free right module. This allows to define partial
derivatives (vector fields).
Corresponding calculi and a differential are said to be a {\it
coordinate
calculi} and a {\it coordinate differential}.
Any coordinate differential defines a commutation rule
$vdx^i=dx^k\c A(v)^i_k$, where $A:v \mapsto A(v)^i_k$ is an algebra
homomorphism $A:R \rt R_{n \times n}$. It is easy to see that for any
homomorphism $R \rt R_{n \times n}$ there exists not more than one
coordinate differential.
We have considered the existence problem of
such a differential. We showed that for a given commutation rule
$vdx^i=dx^k\c A(v)^i_k$ a free algebra generated by the variables
$x^1,\ld ,x^n$ has a related coordinate differential.
We are defining an
{\it optimal} algebra with respect to a fixed commutation rule. In the
homogeneous case this algebra is characterized as the unique algebra
which has no nonzero $A$-invariant subspaces with zero differentials.

In the section II we shall recall the general formalism from ref. [1].
In the section III
we shall consider a number of examples of optimal algebras for
given commutation rules.
The last section is devoted to classification of
commutation rules in two variables which determine commutative
optimal algebra.
\medskip

{\bf 2. Coordinate differentials and optimal algebras. }\
In this section we mainly review our results from [1].
The proofs of all
theorems in this section are in ref. [1].

Let $F$ be a field.
Throughout this paper, an algebra means an unital associative
$F$-algebra generated by a $F$-space $V$
with a basis $\1,\ld ,x^n$.
A presentation of an algebra $R$ is an epimorphism
$\pi :F<\x^1,\ld , \x^n>\rt R$, where $F<\x^1,\ld ,\x^n>$ is the free
associative unital algebra generated by the variables $\x^1,\ld
,\x^n$.
Let $I_R=ker\pi$, then $R\cong F<\x^1,\ld , \x^n>/ I_R$, $I_R$
is an (twosided) ideal of relations in $R$ and $\pi(\x^i)=x^i$
for $i=1,\ld , n$.

A differential from an algebra $R$
to a bimodule $M$ is a linear map satisfying the Leibniz rule
$$
d(uv)=(du)v + udv.
$$\medskip

DEFINITION 2.1 \
A differential $d:R\rightarrow M$ is said to be {\it coordinate} if
bimodule
$M$ is a free right $R$-module freely generated
by $dx^1,\ld , dx^n$.\hfill\break\medskip

NOTE: Coordinate calculi were studied in ref. [1] under the name
{\it free calculi}, due to the freeness condition for the bimodule
$M$ of one forms. More general calculi, {\it calculi with partial
derivatives} and their extensions to higher order (quantum de Rham
complexes)
have been recently considered by the same authors in
[4].\hfill\break\medskip

If $d$ is a coordinate differential, then a linear maps $D_k:R\rt R$,
partial
derivatives, are uniquely defined by the formula,
$$
d\,v=dx^k\c D_k(v)
\eqno          (2.1)
$$
Then
$$
D_k(x^i)= \delta^i_k ,
\eqno     (2.2)
$$
where $\delta^i_k$ is the Kronecker delta.\medskip

PROPOSITION 2.2 \ {\it
A linear map $A_d:R\rt R_{n\times n}$ from an algebra $R$ into the
algebra of $n$ by $n$ matrices over $R$ given by the formula
$$
A_d(v)^i_k\equiv D_k(vx^i)- D_k(v)x^i  \eqno  (2.3)
$$
is an algebra homomorphism i.e.
$$
A_d(uv)^i_k =A_d(u)^l_k A_d(v)^i_l  \eqno  (2.4)
$$
The partial derivatives  $D_k$ and a homomorphism $A_d$ are connected
by the relation
$$
D_k(uv)=D_k(u)v + A_d(u)_k^iD_i(v)   \eqno    (2.5)
$$ }\medskip

The inverse statement also holds.\medskip

PROPOSITION 2.3 \ {\it
Let $R$ be an  algebra generated by  $x^1,\ld ,x^n$
and $A:R \rt R_{n \times n}$ be an algebra homomorphism. If $D_k:R\rt
R,
\ k=1,\ld ,n$ are linear maps such that
$$
D_k(x^i)=\delta^i_k  \eqno   (2.6)
$$
$$
D_k(uv)=D_k(u)v+A(u)_k^iD_i(v), \eqno   (2.7)
$$
then the map $\Delta:v\mapsto dx^k\c D_k(v)$ is a coordinate
differential,
where $M_\Delta=\sum dx^i\c R$ is a free right module with the
left module structure defined by commutation rule $vdx^i=dx^kA(v)^i_k$
i.e. $A_\Delta =A$.}\hfill\break\medskip

Let $M^*$ be a free left $R$-module freely generated by the partial
derivatives $D_i$, i.e. $Y\in M^*$ iff $Y=Y^iD_i$ with $Y^i\in R$.
Define a right $R$-module structure on $M^*$ by the transpose
commutation rule
$$
Y.v\equiv Y^i(D_i .v)= Y^iA(v)_i^kD_k
$$
Thus we have defined a bimodule $M^*$ of vector fields as a dual
to a bimodule $M$ of differential forms together with a
pairing $<Y, \omega>\equiv Y^i\omega_i \in R$, where $\omega =
dx^i\c\omega_i
\in M$. Then $<Y.f, \omega>=<Y, f\omega>$.
A vector field $Y\in M^*$ can be characterized as a linear map
(endomorphism) $Y:R\rt R$ satisfying the twisted Lebniz's rule vis.
$$
Y(uv)= Y(u)v + (Y.u)(v) \eqno (2.8)
$$
where $Y(u)=Y^iD_i(u)$ and $Y.u$ is a previously defined
multiplication
from the right. It generalizes the formulae (2.5) to arbitrary vector
field $Y$.\hfill\break
Both definitions  differential forms and vector fields essentially
depend
on the generating space $V={\rm lin}(x^1,\ld ,x^n)$ in the following
sense.
Let $z^k= \al^k_i x^i$ be another base in $F$-space $V$
with a matrix $(\al^k_i)\in GL(n, F)$. Then the basis of a bimodule
$M,$
a differentials
$dx^i$ and basic vector fields $D_i= \partial / \partial x^i$
undergo correspondingly covariant and contravariant transformation
law, i.e.
$$
dz^k = \al^k_i dx^i,\ \ \ \
\partial /\partial z^k = \be_k^i \partial / \partial x^i
$$
where $\be^k_l\al^l_i=\de^k_i$.\hfill\break
A  question concerning Proposition 2.3 arises here. If a
homomorphism $A$ is given, then formula (2.7) allows one to calculate
partial derivatives of a product in terms of its factors. That fact
and
formula (2.6) show that for a given $A$ there exists not more then
one $D$
satisfying formulas (2.6) and (2.7). It is not clear yet whether or
not there
exists at least one $D$ of such a type. Thus, our first
task is to describe these homomorphisms of $A$ for which there exist
coordinate
differentials with $A_d=A$.\hfill\break\medskip

THEOREM 2.4 \ {\it
Let $R=F<x^1,\ld ,x^n>$ be a free  algebra generated by
$x^1,\ld ,x^n$ and $A^1,\ld ,A^n$ be any set of $n\times n$ matrices
over $R$. There exists the unique coordinate differential $d$ such
that
$A_d(x^i)=A^i$.}\hfill\break\medskip

Let $R$ be an $F$-algebra defined by a set of
generators \ \ \ $x^1,\ld ,x^n$\ \ and a set of relations \
$f_m(x^1,\ld ,x^n)\ =\ 0,\ \ m\in \M$ \ i.e.\ $R=\R\,/\, \I$,\ \
where \
$\R=F<\x^1,\ld ,\x^n>$ is a free algebra  and $\I$ is its ideal
generated by elements $f_m(\x^1,\ld ,\x^n),\ \ m\in \M$.

Let $\pi$ be an algebra projection $\R \rt R$ such that
$\pi(\x^i)=x^i$, $\pi(1)=1$. Since $R_{n\times n}= R\ot F_{n\times
n}$,
$\pi$ defines an epimorphism $\pih:\R_{n\times n}\rt R_{n\times n}$
by
the formula $\pih=\pi \ot id.$

If $A:R\rt R_{n\times n}$ is an algebra map, then we have
the following diagram of algebra homomorphisms
$$
\matrix{\R &{\buildrel\A\over\longrightarrow} & \R\ot F_{n\times n}&=
&\R_{n\times n}\cr
\pi \downarrow &   & \downarrow \pi \ot {\rm id} & = & \downarrow
\pih\cr
R &{\buildrel A\over\longrightarrow} & R\ot F_{n\times n} & = &
R_{n\times n}}\eqno(2.9)$$

Let choose for each generator $x^i$ an arbitrary element $\A^i\in \R$
such that $\pih(\A^i)=A^i$. Then the
map $\x^i\mapsto  \A^i$ can be extended to an algebra homomorphism
$\A:\R
\rt \R_{n\times n}$ (recall that $\x^1,\ld ,\x^n$ are free
variables).
This homomorphism completes a diagram (2.9) to a commutative diagram.
For each relation
$f_m(\x^1,\ld ,\x^n)$ we have:
$$
\pih(\A(f_m)) = A(\pi(f_m))=0. \eqno (2.10)
$$
Furthermore,
$$
ker\pih=ker(\pi\ot id)=ker\pi\ot F_{n\times n}=\I_{n\times n},
$$
and finally,
$$
\A(f_m) \in ker\pih=\I_{n\times n}.
$$
Theorem 2.4 claims that for a homomorphism $\A:\R\rt \R_{n\times n}$
there exists a unique coordinate differential $\d$ of a free algebra
$\R$.\hfill\break\medskip

DEFINITION 2.5 \
The differential $\d$ is called a {\it cover} differential with
respect to the homomorphism $A:R\rt R_{n\times
n}$.\hfill\break\medskip

We have proved:\medskip

THEOREM 2.6 \ {\it
For any homomorphism $A:R\rt R_{n\times n}$ there exists a cover
differential $\d$ of a free algebra $\R$.}\hfill \break\medskip

PROPOSITION 2.7 \ {\it
Let an algebra $R$ be  generated by $x^1,\ld ,x^n$ subject to the
relations $\{f_m=0,\,m\in \M \}.$
Let $\hat{D}_k$ be a partial derivatives of the cover differential
$\d$,
and $\I$ be an ideal generated by $\{f_m,\ m\in \M \}$. Then an
algebra
$R$ possesses a coordinate differential with
respect to a homomorphism $A:R\rt R_{n\times n}$ if and only if}
$$\hat{D}_k(f_m(\x^1,\ld ,\x^n))\in \I$$\medskip

COROLLARY 2.8 \ {\it
Let an algebra $R$ be  generated by $x^1,\ld ,x^n$ subject to
the set of homogeneous relations  $\{f_m=0, \ m\in \M\}$ of the same
degree.
If $A:R\rt R_{n\times n}$ acts linearly on generators $A(x^j)^i_k=
\al^{i j}_{kl}x^l$, then for a pair ${R, A}$ there exists a coordinate
differential if and only if for all $m\in \M$}
$$
\d f_m=0\ \ .  \eqno  (2.11)$$\hfill\break\medskip

DEFINITION 2.9 \
An ideal $I\neq \R$ of a free algebra $\R=F<\x^1,\ld ,\x^n>$ is said
to
be {\it consistent} with a homomorphism $A:\R\rt \R_{n\times n}$ if a
factor algebra $\R/I$ has a coordinate differential satisfying the
commutation
rules
$$
x^jdx^i=dx^k\c A(x^j)^i_k\ \ .  $$\hfill\break \medskip

If an ideal $I$ is $A$-consistent, then Proposition 2.2 defines a
homomorphism
$A:r\mapsto A^i_k(r)$ from the factor algebra into the matrix algebra
over it. Thanks to Proposition 2.7 it follows that $I$ is
$A$-{\it consistent}
iff $A^i_k(I)\subseteq I$,
and if $D_k(I)\subseteq I$ for
any of partial derivatives $D_k$ defined by a differential $d$
corresponding to $A$ (see Theorem 2.4).\hfill\break
These two condition generalize
Wess and Zumino's quadratic and a linear consistency conditions for
quadratic algebras [7].\hfill\break

We have proved that free algebra $\R$ admits a coordinate
differential for
arbitrary commutation rules. In order to define a homomorphism
$A:\R\rt
\R_{n\times n}$ it is enough to set its value on generators vis.
$$
A^i_k(x^m)=\al^{mi}_k +\al^{mi}_{kj_1}x^{j_1}+
\al^{mi}_{kj_1j_2}x^{j_1}x^{j_2}+ \ld
$$
where $\{\al^{mi}_{kj_1,\ld ,j_r} \}$ are arbitrary tensor
coefficients.
If a homomorphism $A$ preserves a degree, then it must act linearly on
generators $A^j_k(\x^i)=\al^{ij}_{kl}\x^l$. Therefore, a homomorphism
$A$ is defined by a 2-covariant and 2-contravariant tensor
$A=\al^{ij}_{kl}$. This is a homogeneous case.
The case $A^j_k(x^i)= \al^{ij}_k+
\al^{ij}_{kl}x^l$ has been considered by Dimakis and M\"uller--
Hoissen [2].
In this case the partial derivatives do not increase the degree.
\hfill\break

For any homomorphism $A$ there exists the largest $A$-consistent
ideal
$I(A)$ contained in the ideal $\bar{R}$ of polynomials with zero
constant
terms ($\bar{R} \dot{+} F\c 1=\R$) -- the sum of all consistent
ideals of
such a type. \medskip

NOTE: \  An ideal $I(A)$ need not to be the only
maximal $A$-consistent ideal in $\R$.\medskip

THE MAIN DEFINITION. \
The factor algebra $R_A=\R / I(A)$ is said to be an
{\it optimal} algebra for the commutation
rule $x^idx^j=dx^k\c A(x^i)^j_k$. A pair: the optimal algebra $R_A$,
and the differential $d$
corresponding to the commutation rule $A$
will be called an {\it optimum calculus}.\medskip

We will describe an ideal $I(A)$ in the homogeneous case.\medskip

THEOREM 2.10 \ {\it
For each 2-covariant 2-contravariant tensor $A=\al^{ij}_{kl}$
the ideal $I(A)$ can be constructed by induction as the homogeneous
space
$I(A)=I_1(A)+I_2(A)+I_3(A)+\cdots$ in the following way:

   \item{(i)}  $I_1(A)=0$
   \item{(ii)}  Assume that $I_{s-1}(A)$ has been defined and $U_s$
be a space
   of all polynomials $m$ of degree $s$ such that $D_k(m)\in I_{s-
1}(A)$
  for all $k.\ 1\leq k\leq n$. Then $I_s(A)$ is the largest $A$-
invariant
    subspace of $U_s$.

The ideal $I(A)$ is a maximal $A$--consistent ideal in $\R$.}\medskip

COROLLARY 2.11 \ {\it
The Theorem 2.10 shows in particular that
if a homogeneous element is such that all elements of the invariant
subspace generated by it have all partial derivatives equal to zero,
then
that element vanishes in the optimal algebra.}\medskip

{\bf 3. Examples of optimum calculi. }\
In this section we shall consider a number of commutation rules $A$
and the
corresponding optimal algebras $R_A$.\medskip

EXAMPLE 3.1. \ {\it
Consider a diagonal commutation rule, $x^jdx^i=dx^i\c
q^{ij}x^j$, with $q^{ij}q^{ji}=1, i\neq j$. If no one
of the coefficients $q^{ii}$ is a root of a polynomial of the type
$\la^{[m]}\doteq \la^{m-1}+\la^{m-2}+\cdots +1$, then the optimal
algebra is $R_A= F<x^1,\ld ,x^n>/\{q^{ij}x^ix^j=x^jx^i, \ i<j
\}$.\hfill\break
If $(q^{ii})^{[m_i]}=0,\ 1\leq i\leq s$ with minimal $m_i$ then
\hfill\break
$R_A=F<x^1,\ld ,x^n>/\{q^{ij}x^ix^j=x^jx^i,\ i<j,\ (x^i)^{m_i}=0,\
1\leq i\leq s \}$.}

{\it Proof. }\
First of all we have to note that the elements $q^{lj}x^lx^j-
x^jx^l$ are zero in the optimal algebra. By Theorem 2.10 it is
sufficient
to show that all partial derivatives (for cover differential) of
elements
of the invariant space generated by this forms are equal to zero. In
our
example $A(x^j)^i_k=\de^i_k q^{ij}x^j$ and therefore
$$
A(q^{lj}x^lx^j-x^jx^l)^i_k=q^{lj}A(x^l)^i_sA(x^j)^s_k -
A(x^j)^i_sA(x^l)^s_k =
$$
$$
=q^{lj}\de^i_sq^{il}x^l\c \de^s_kq^{sj}x^j -
\de^i_sq^{ij}x^j\c \de^s_kq^{sl}x^l =
$$
$$
=\de^i_kq^{kl}q^{kj}[q^{lj}x^lx^j - x^jx^l] \eqno (3.1)
$$
So it is enough to check that the partial derivatives of the relations
are zero:
$$
D_k(q^{lj}x^lx^j-x^jx^l)=q^{lj}[D_k(x^l)x^j+A(x^l)^i_kD_i(x^j)]-
$$
$$
-[D_k(x^j)x^l-A(x^j)^i_kD_i(x^l)]=0$$
Formula (3.1) and Corollary 2.7 show that the factor algebra
$S=F<x^1,\ld,x^n>/\{q^{ij}x^ix^j=x^jx^i,\ \ i<j\}$ has coordinate
differential
with our commutation rules. For this algebra to be optimal by Theorem
2.10
it is enough to see that any homogeneous element of positive degree
with
zero all the partial derivatives is equal to zero in this
algebra.\hfill\break
We have
$$
D_k[(x^j)^m]=\de^j_k(x^j)^{m-1} + A(x^j)^i_kD_i[(x^j)^{m-1}]=
$$
$$
=\de^j_k(x^j)^{m-1} + \de^i_kq^{ij}x^jD_i[(x^j)^{m-1}]
$$
and by an easy induction
$$
D_k[(x^j)^m]=\de^j_k[1+q^{jj}+(q^{jj})^2+\ld +(q^{jj})^{m-1}](x^j)^{m-
1}\ .
\eqno (3.2)
$$
An arbitrary element of the algebra $S$ has a unique presentation
of the form
$$
f=\sum \al_i(x^1)^{i_1}(x^2)^{i_2}\ld (x^n)^{i_n}\ ,
$$
where, $i\equiv i_1i_2\ld i_n$ is a multi index. Thus by formula (3.2)
$$
D_kf=\sum \al_i(x^1)^{i_1}(x^2)^{i_2}\ld D_k[(x^k)^{i_k}]\ld
(x^n)^{i_n}\ =
$$
$$
=\sum (q^{kk})^{[i_k-1]} \al_i(x^1)^{i_1}(x^2)^{i_2}\ld
(x^k)^{i_k-1}\ld (x^n)^{i_n}\ . \eqno (3.3)
$$
If no one of the elements $(q^{kk})^{[m]}$ is zero, then the right
hand side
components of the formula (3.3) are equal to zero for all $k$ if and
only if
$\al_i=0$. So $f=0$ and $S$ is the optimal algebra.

If $(q^{ii})^{[m_i]}=0,\ 1\leq i\leq s$, then by (3.2) we have
$D_k[(x^i)^{m_i}]=0$ and also $D_k[A\{(x^j)^{m_j}\}]=0$ because of
$A(x^j)^i_k=\de^i_kq^{ij}x^j$ implies that $A\{(x^j)^{m_j}\}$ has the
form $\Lambda\c (x^j)^{m_j}$, where $\Lambda$ is a matrix with
coefficients
from the base field. Therefore $(x^j)^{m_j}=0$ in the optimal
algebra.\hfill\break
Each element of the algebra
$$
\bar{S}=F<x^1,\ld ,x^n>/\,\{q^{ij}x^ix^j=x^jx^i,\ i<j,\ (x^i)^{m_i}
=0,\ 1\leq i\leq s\}
$$
has unique presentation of the form (3.1), where
$i_1<m_1, i_2<m_2,\ld ,i_s<m_s$. Now formulae (3.2) and (3.3), which
are
still valid in $\bar{S}$, imply that if all the partial derivatives of
an element are zero in $\bar{S}$ then this element is zero in
$\bar{S}$.
\hfill Q.E.D.\hfill\break\medskip

EXAMPLE 3.2. \ {\it
Let $A=0$ i.e. $x^idx^j=0$. Then $\d$ is a homomorphism of right
modules and the optimal algebra is free $R_A=\R$.}
\hfill\break

{Proof: }\ Indeed, if $u=x^1u_1+\ld +x^nu_n$ then $D_k(u)=u_k$ and
any element
with zero partial derivatives is zero. Thus by Theorem 2.10 the ideal
$\I_A$ contains no nozero elements.\hfill Q.E.D.\hfill\break\medskip

Note that the calculus with this commutation rule in fact was defined
in Fox's paper [3]. This calculus is essential tool in series of
Fox's papers on group theory. It  was defined for the free group
algebra
by the following formula for derivations:
$$
D(uv)=D(u)v\,+\,u^o D(v)
$$
where, $u^o$ is the sum of all coefficient of the element $u$ from
group
algebra of free group $G$ freely generated by the elements $x^1,\ld
,x^n$.

If partial derivatives $D_k$ for a differential $d$ obey this
conditions
then the commutation rule has the form $x^idx^j=dx^j$. If we change
the
variables $X^i=x^i-1$ than we will obtain the calculus with zero
homomorphism $A$.\medskip

EXAMPLE 3.3. \ {\it
Let $x^idx^j=-dx^i\c x^j$. Then the optimal algebra is the
smallest possible algebra generated by the space $V$ i.e.
$R_A=F<x^1,\ld
,x^n>/ \{x^ix^j=0 \}=F+V$.}

{\it Proof: }\ In this case $\d (x^ix^j)=\d x^i\c x^j+x^i\d x^j=0$
for cover
differential. Evidently the space
of all quadratic forms is $A$--invariant. By Theorem 2.10  in
the optimal algebra $x^ix^j=0$.\hfill Q.E.D.\hfill\break\medskip

EXAMPLE 3.4. \ {\it
Let $x^1dx^1= dx^1\c (\al_2x^2+\cdots +\al_nx^n)$ and
$x^idx^j=-dx^i\c x^j$ if $i\neq 1\ or \ j\neq 1$. Then the optimal
algebra
is almost isomorphic to the ring of polynomials in one variable. More
precisely,
$R_A=F<x^1,\ld ,x^n>/ \{x^ix^j=0, \ unless \ \,i=j=1 \}$.}

{\it Proof. }\
In this case $A^1_k(x^1)=\de^1_k(\al_2x^2+\ld +\al_nx^n)$ and
$A^j_k(x^i)= -\de^i_kx^j$ if $i\neq 1$ or $j\neq 1$.
Let $I$ be an ideal generated by all products $x^ix^j$ but $x^1x^1$
and let $S$ be a factor algebra $\R/I$. We have
$$
\D_k(x^ix^j)=\de^i_kx^j+A^s_k(x^i) \de^j_s = \de^i_kx^j+A^j_k(x^i)=0
$$
if either $i\neq 1\ \ or\ \ j\neq 1$. We also have
$$
A^m_s(x^ix^j)=A^m_k(x^i)A^k_s(x^j)
$$
Both left and right factors contain an addendum of the form $\al x^1$
only if both $m=k=1$ and $i=k,\ \ j=s$ are valid and both $m=i=1$ and
$j=k=1$ are false, which is impossible. Thus the product belong to the
ideal $I$.
By the Proposition 2.7 $I$ is a consistent ideal.\hfill\break
Any element of the factor algebra $S$ has unique presentation of the
form
$$
f=\ga_kx^k + \be_2(x^1)^2+\ld +\be_N(x^1)^N
$$
For $k\neq 1$,  we have
$$D_k[(x^1)^N]= A^s_k(x^1)D_s[(x^1)^{N-1}]= -x^1D_1[(x^1)^{N-1}]
$$
and
$$D_1[(x^1)^N]=(x^1)^{N-1} + A^k_1(x^1)D_k[(x^1)^{N-1}] =
(x^1)^{N-1} + w\c D_1[(x^1)^{N-1}] -
$$
$$- \sum_{k\geq 2}x^k\c D_k[(x^1)^{N-1}]=
\left\{\matrix{(x^1)^{N-1} & \hbox{if $N\geq 3,$}\cr
x^1 + w & \hbox{if $N=2$,}\cr}\right.$$
where $w=\al_2x^2+\ld + \al_nx^n$. Now if $df=0$ in $S$ then
$\ga_k=D_k(f)=0$ for $k\geq 2$ and
$$
D_1f=\ga_1+\be_2x^1+\be_2w+\be_3(x^1)^2+\ld +\be_N(x^1)^{N-1}=0
$$
Therefore $\be_1=\ld =\be_N=0$ and $R_A=S$.
\hfill Q.E.D.\hfill\break\medskip

EXAMPLE 3.5. \ {\it
Let $\mu, \la \in F$.
If $n=2$ and $x^1dx^1=dx^1\c \mu x^2, \ \ x^1dx^2=-dx^1\c x^2, \ \
x^2dx^1=-dx^2x^1, \ \ x^2dx^2=dx^2\c \la x^1$, then the optimal
algebra is
isomorphic to the direct sum of two copies of the polynomial algebra
$R_A=F<x^1, x^2>/ \{x^1x^2=x^2x^1=0 \}$.}\hfill\break

{\it Proof. }\
Let $I$ be the ideal generated by  $x^1x^2$ and $x^2x^1$. We have
$\d (x^1x^2)=\d (x^2x^1)=0$. By definition
of this commutation rule we have
$$
A(x^1)=\pmatrix{\mu x^2 & -x^2 \cr  0 & 0\cr} ;\qquad
A(x^2)=\pmatrix{0 & 0 \cr -x^1 & \la  x^1\cr} ;$$
$$
A(x^1x^2)=\pmatrix{x^2x^1 & -\la x^2x^1 \cr  0 & 0\cr}\equiv 0\quad
(mod\; I);
$$
$$
A(x^2x^1)=\pmatrix{0 & 0 \cr -\mu x^1x^2 & x^1x^2\cr}
\equiv 0 \quad (mod\; I).
$$
Therefore the ideal $I$ is consistent.\hfill\break
In the algebra $S\equiv F<x^1, x^2> / I$ any element has a unique
presentation
of the form
$$
f=\al_1x^1+\al_2(x^1)^2+\ld +\al_n(x^1)^n +\be_1x^2+\be_2(x^2)^2+
\ld +\be_m(x^2)^m
$$
We have $D_2[(x^1)^n]= A^k_2(x^1)\c D_k[(x^1)^{n-1}]=0$ and therefore
$$
D_1[(x^1)^n]=(x^1)^{n-1}+A_1^k(x^1)\c D_k[(x^1)^{n-1}]=
$$
$$
=(x^1)^{n-1}+\mu x^2\c D_1[(x^1)^{n-1}]-x^2\c D_2[(x^1)^{n-1}]
=(x^1)^{n-1} + (\mu x^2)^{n-1}
$$
By this formula
$$
D_1f=\al_1+\al_2x^1+\ld +\al_n(x^1)^{n-1}+\al_2\mu x^2+
\al_3(\mu x^2)^2+\ld +\al_n(\mu x^2)^{n-1}
$$
Therefore $D_1f=0$ only if $\al_1=\al_2=\ld =\al_n=0$. Analogously,
$D_2f=0$ only if $\be_1=\be_2=\ld =\be_m=0$. By Theorem 2.10
in this case the optimal algebra is $S$.\hfill
Q.E.D.\hfill\break\medskip

{\bf 4. Homogeneous commutation rules in two variables
with commutative optimal algebra. }\
In this
section we will describe all homogeneous commutation rules in
two variables with a commutative optimal algebra. In this case the
ideal
$I(A)$ is homogeneous and is defined by Theorem 2.10. Commutativity
of the
optimal algebra is equivalent to $x^1x^2-x^2x^1 \in I_2$, where $I_2$
is the second homogeneous component of $I(A)$. We will call the
commutation
rule (and the corresponding optimal algebra) {\it regular} if the
space
$I_2$ is one dimensional, i.e. if it is generated by commutator. In
the
opposite case we will call the commutation rule and the optimal
algebra
{\it irregular}. For instance the optimum calculi in the Examples 3--5
are irregular. Evidently if the optimal algebra is isomorphic to the
algebra of polynomials in two variables then the commutation rule is
regular (but not vice versa).\hfill\break\medskip

THEOREM 4.1 \ {\it
Let $u,\ \ v,\ \ w,\ \ v_1\in V=\ {\rm lin}(\1, \2)$ and $\la, \
\mu\in F$.
A homogeneous commutation rule with
regular commutative optimal algebra belongs (up to renaming of
variables
$\1\leftrightarrow\2$) to one of the following four classes:}

 \item{(I)}
$$\eqalign{x^1dx^1 & = dx^1\c u\,+\,dx^2\c v, \cr
x^1dx^2 & =dx^1\c w\,+\,dx^2\c (\la v+x^1),\cr
x^2dx^1 & =dx^1\c(w+x^2)\,+\,dx^2\c(\la v),\cr
x^2dx^2 & =dx^1\c (\la w)\,+\,dx^2\c (\la^2v-\la u+w+\la
x^1+x^2);\cr}$$

  \item{(II)}
$$\eqalign{x^1dx^1 & =dx^1\c (\1+\mu v+v_1)\,+\,dx^2\c v, \cr
x^1dx^2 & =dx^1\c \la v\,+\,dx^2\c (v_1+x^1),\cr
x^2dx^1 & =dx^1\c(\la v+x^2)\,+\,dx^2\c v_1,\cr
x^2dx^2 & =dx^1\c \la v_1\,+\,dx^2\c (\la v-\mu v_1+x^2);\cr}$$

  \item{(III)}
$$\eqalign{x^1dx^1 & =dx^1\c u,\cr
x^1dx^2 & =dx^2\c x^1,\cr
x^2dx^1 & =dx^1\c x^2,\cr
x^2dx^2 & =dx^2\c v;\cr}$$

 \item{(IV)}
$$\eqalign{x^1dx^1 & =dx^1\c u,\cr
x^1dx^2 & =dx^2\c u,\cr
x^2dx^1 & =dx^1\c x^2\,+\,dx^2\c (u-\1),\cr
x^2dx^2 & =dx^1\c w\,+\,dx^2\c v.\cr}$$

{\it Proof. }\
First of all we will prove that each commutation rule has a
commutative
optimal algebra.\hfill\break
Let a commutation rule is set by formulae I. It is enough to prove
that
the ideal $I$ generated by the commutator $x^1x^2-x^2x^1$ is
consistent.
We have
$$A(x^1)=
\pmatrix{u & w \cr v & \la v+x^1};$$
$$A(x^2)=\pmatrix{w+x^2 & \la w \cr \la v & \la^2v-\la u+w+\la
x^1+x^2}
\eqno (4.1)$$
Therefore $A(x^2)=\la A(x^1)+(w+x^2-\la u)E$, where $E$ is
a unit matrix and evidently
$$
A([x^1, x^2])=[A(x^1), A(x^2)]=[A(x^1), \la A(x^1)+(w+x^2-\la u)E]=0
$$
in the factor algebra $F<x^1, x^2>/\{x^1x^2=x^2x^1\}$, so
$A^i_kI\subset
I.$
For partial derivatives we have
$$
D_1(x^1x^2-x^2x^1)=x^2+A^k_1(\1)D_k(\2)-A^k_1(\2)D_k(\1)=
$$
$$
=\2 +A^2_1(\1)-A^1_1(\2)=\2 +w-(w+\2)=0,
$$
$$
D_2(\1 \2 -\2 \1)=A^k_2(\1)D_k(\2)-\1 - A^k_2(\2)D_k(\1)=
$$
$$
=A_2^2(\1)-\1 -A^1_2(\2)=\la v+\1 -\1 -\la v=0
$$
and all the commutation rules from series (I) have commutative optimal
algebra.\hfill\break

Let a commutation rule is set  by formula (II). We have
$$A(\1)=\pmatrix{\1 +\mu v+v_1 & \la v \cr v & \1 +v_1},$$
$$A(\2)=\pmatrix{\2 +\la v & \la v_1 \cr v_1 & \2 +\la v-\mu v_1},
\eqno (4.2)$$
$$A([\1, \2])=[A(\1), A(\2)]=0.$$
Therefore $A^i_kI\subset I.$ Moreover
$$D_1(\1 \2 -\2 \1)=\2 +A^2_1(\1)-A^1_1(\2)=0,$$
$$D_2(\1 \2 -\2 \1)=A^2_2(\1)-\1 -A^1_2(\2)=0,$$
and all the commutation rules given by the formulae (II) have
commutative
optimal algebra.

Let a commutation rule is set by the formulae (III). We have
$$
A(\1)= \pmatrix{u & 0 \cr 0 & \1},\quad
A(\2)= \pmatrix{\2 & 0 \cr 0 & v}.\eqno (4.3)$$
Because $A(\1)A(\2)=A(\2)A(\1)$ then $A^i_kI\subset I.$ Moreover
$$D_1(\1 \2 -\2 \1)=\2 +A^2_1(\1)-A^1_1(\2)=0,$$
$$D_2(\1 \2 -\2 \1)=A^2_2(\1)-\1 -A^1_2(\2)=0,$$
and all the commutation rules given by formulae (III) have commutative
optimal algebra.\hfill\break
In the case (IV) we have
$$A(\1)= \pmatrix{u & 0 \cr 0 & u},\quad
A(\2)= \pmatrix{\2 & w \cr u-\1 & v}.\eqno (4.4)$$
It is evident that $A^i_kI\subset I.$ Moreover
$$D_1(\1 \2 -\2 \1)=\2 +A^2_1(\1)-A^1_1(\2)=0,$$
$$D_2(\1 \2-\2 \1)=A^2_2(\1)-\1 -A^1_2(\2)=0,$$
and all of the commutation rules in the theorem have commutative
optimal
algebra.\hfill\break
Conversely, let a commutation rule $x^idx^j=dx^kA^{ij}_k$, where
$A^{ij}_k
=A^j_k(x^i)$, has commutative optimal algebra $S$. In this case
$$
0=D_k(x^ix^j-x^jx^i)=\de^i_kx^j+A^{ij}_k-\de^j_kx^i-A^{ji}_k
$$
As a consequence we obtain necessary conditions (and they hold for
$n>2$):
$$
A^{ij}_k=A^{ji}_k\ \ \ \ \hbox{if $k\neq i$ and $k\neq j, $}\eqno
(4.5)
$$
$$
A^{ij}_j=x^i+A^{ji}_j\ \ \ \ \hbox{if $j\neq i$ and $k=j.$}\eqno (4.6)
$$
If $n=2$ then these conditions are reduced to the following two
$$
A^{12}_2=\1 +A^{21}_2 ;\ \ \ \ \ \ \ A^{21}_1=\2 +A^{12}_1 \eqno (4.7)
$$
Also we have $A(\1 \2 -\2 \1)=0$ i.e. the matrices $A^i=A^{is}_k$
and $A^j=A^{js}_k$ commutes in the ring of matrices over $S$:
$$
A^{is}_kA^{jk}_m=A^{js}_kA^{ik}_m \eqno (4.8)
$$
Let  consider these equalities in detail. All $A^{is}_k$
have degree one i.e. they are in the space $V$. So the relations
(4.8) have
degree  two and therefore  have to belong to the second
homogeneous component $I_2$ of the ideal $I(A)$. Let $S_2$ be a factor
algebra $\R/ I_2$. This is a commutative algebra and relations the
(4.8) are
valid on it. As $S$ is regular, the space $I_2$ is generated by
commutator
$\1 \2 -\2 \1$ and the algebra $S_2=F[\1, \2]$ is the algebra of
polynomials
in two variables.\hfill\break
If one of the matrices $A(\1), \, A(\2)$ is scalar then (if necessary
by
renaming  variables) we can suppose that
$ A(\1)= \pmatrix{u & 0 \cr 0 & u\cr}$ and the relations (4.7)
show that $A(\2)^1_1= \2, \ \ A(\2)^1_2=u-\1$ and the commutation rule
belongs to the series (IV). If both matrices $ A(\1), A(\2)$ are
diagonal
than relations (4.7) immediately imply
$A(\1)= \pmatrix{u & 0 \cr 0 & \1\cr},\quad A(\2)=
\pmatrix{\2 & 0 \cr 0 & v\cr}$ and the commutation rule
belongs to the series (III). So under the following consideration we
can
suppose that no one of the matrices $A(\1), A(\2)$ is scalar and one
of them is not diagonal matrix.\hfill\break
In the algebra of $2 \times 2$ matrices over the field of rational
functions
$K=F(\1, \2)$ the dimension over the field $F$ of a centralizer of any
non scalar matrix is equal to $2$. It means that the centralizer of
the
matrix $A(\1)$ is generated (over $K$) by two matrices $A(\1)$ and
$E$. It implies that we have the relation $A(\2)=g\c A(\1)+f\c E$ with
$f, g \in K$. It implies $A^{22}_1=gA^{12}_1,\ \ A^{21}_2=gA^{11}_2$
and therefore $A^{21}_2\c A^{12}_1=A^{22}_1A^{11}_2$. By definition
all
the coefficients in the matrices are linear combinations of the
variables
so either $A^{21}_2=\la A^{11}_2, A^{22}_1=\la A^{12}_1$ or
$A^{12}_1=\la A^{11}_2, \ \ A^{22}_1=\la A^{21}_2$, where $\la \in F$
or $\la=\infty$. The last case $\la =\infty$ means respectively
$A^{11}_2=A^{12}_1=0$ or $A^{11}_2=A^{21}_2=0$. These two cases
reduce to
the cases $\la =0$ by changing the variables $\1 \leftrightarrow
\2$.\hfill\break
If $A^{21}_2=\la A^{11}_2, \ \ A^{22}_1=\la A^{12}-1$, then a
denotation
$A^{11}_1=u, \ A^{12}_1=w$ implies $A^{21}_2=\la v, \ A^{22}_1=\la w$
and therefore $f=\2 +w-\la u$. So $A^{22}_2=\la A^{12}_2+f=\la \1 +
\la^2v+\2+w-\la u$. It means that the matrices $A(\1), \, A(\2)$ have
the form~(~I~).\hfill\break
If $A^{12}_1=\la A^{11}_2, \ A^{22}_1=\la A^{21}_2$ then it is
convenient
to denote $A^{11}_2=v, \ A^{21}_2=v_1$. In this case $g={v_1\over v}$
and $f=A^{21}_1-{v_1\over v}A^{11}_1=A^{22}_2-{v_1\over v}A^{12}_2$
or $v_1(A^{11}_1-A^{12}_2)=v(A^{21}_1-A^{22}_2)$. All the factors in
the
last relation are linear combinations of the variables, therefore
$A^{11}_1-A^{12}_2=\mu v,\ \ A^{21}_1-A^{22}_2=\mu v_1$ where $\mu
\in F$
($\mu=\infty$ as $v \neq 0$ or $v_1 \neq 0$ because one of the
matrices
$A^1,\ \ A^2$ is not diagonal). Now relations (4.7) have the form
$A^{12}_2=x^1+v_1$, $A^{21}_1=x^2+\la v$ which implies
$A^{11}_1=x^1+v_1 +\mu v$, and $A^{22}_2=x^2+\la v-\mu v_1$. This
means
that we are in the situation (II). The theorem is proved.\hfill Q.E.D.
\hfill\break\medskip

NOTE: \ It is an open problem to determine the optimal
algebra
for the commutation rules described above. We do not claim
that the optimal algebra is a polynomial algebra in two variables
$S_2=
F<\1, \2>/<I_2>\equiv F[\1, \2]$.\hfill\break\medskip

{\bf 5. Acknowledgments.}\ The authors are greatly indebted to
Zbigniew
Ozie\-wicz for his active interest in the publication of this paper
and for
many stimulating conversations.\hfill\break

\ninerm
\centerline{\bf References}
\item{[1]}
Borowiec A., Kharchenko V. K. and Oziewicz Z., (1993), {\it
On free differentials on associative algebras},
in {\it Non Associative Algebras and Its Applications},
series {\it Mathematics and Applications},
ed. S. Gonz\'alez, Kluwer Academic
Publishers, Dordrecht, 1994, 43-56 (hep-th/9312023)\hfill\break
Borowiec A., Kharchenko V. K., (1994), {\it Coordinate calculi on
associative algebras}, in {\it Quantum Groups: Formalism and
Applications}, Ed. J. Lukierski et al, PWN, Warsaw 1994, 231-242
(hep-th/9501051)

\item{[2]}
Dimakis A. and M\"uller-Hoissen F., (1992), {\it
Quantum mechanics as non-commutative symplectic geometry}, J. Phys.
A: Math.
Gen. {\bf 25} 5625-5648\hfill\break
Dimakis A., M\"uller-Hoissen F. and Striker T., (1993),
{\it Non-commu\-ta\-tive differential calculus and lattice gauge
theory},
J. Phys. A: Math. Gen. {\bf 26} 1927-1949

\item{[3]}
Fox Ralph H. (1953), {\it Free differential calculus. I. Derivation
in the
free group ring}, Annals of Mathematics {\bf 57} (3) 547-
560\hfill\break
Fox Ralph H. (1954), {\it Free differential calculus II.},
Annals of Mathematics {\bf 59} (2) 196-210\hfill\break
Fox Ralph H. (1960), {\it Free differential calculus V.},
Annals of Mathematics {\bf 71} (3) 407-446

\item{[4]}
Kharchenko V. K. and Borowiec A., (1995), {\it Algebraic approach to
calculi with partial derivatives}, to be published in Siberian
Advances
in Mathematics {\bf 5} (2) 1-28

\item{[5]}
Pusz W., (1989), {\it Twisted canonical anticommutation
relations},
Reports on Mathematical Physics {\bf 27} 349-360

\item{[6]}
Pusz W. and  Woronowicz S. L., (1989), {\it Twisted second
quantization}, Reports on Mathematical Physics {\bf 27} 231-257

\item{[7]}
Wess J. and Zumino B., (1990), {\it Covariant differential calculus
on the quantum hyperplane},
Nuclear Physics {\bf 18B} 303-312,\hfill\break Proceedings,
Supplements. Volume in
honor of R. Stora.

\item{[8]}
Woronowicz S. L., (1989), {\it Differential calculus
on compact matrix pseudogroups} (quantum groups), Comm. Math.
Phys. {\bf 122} 125-170.

\item{}\ \ \ \ \ \ \ \ \

\bye